# WEAK CONVERGENCE OF METROPOLIS ALGORITHMS FOR NON-I.I.D. TARGET DISTRIBUTIONS[1]

By Mylène Bédard

*University of Warwick*

In this paper, we shall optimize the efficiency of Metropolis algorithms for multidimensional target distributions with scaling terms possibly depending on the dimension. We propose a method for determining the appropriate form for the scaling of the proposal distribution as a function of the dimension, which leads to the proof of an asymptotic diffusion theorem. We show that when there does not exist any component with a scaling term significantly smaller than the others, the asymptotically optimal acceptance rate is the well-known 0.234.

**1. Introduction.** Metropolis algorithms [8, 9] provide a method for sampling from highly complex probability distributions. The ease of implementation and wide applicability of these algorithms have given them their popularity and they are now frequently used by practitioners at all levels in various fields of application. However, their convergence can sometimes be lengthy, which suggests the need for an optimization of their performance. Because the efficiency of Metropolis algorithms depends crucially on the scaling of the proposal density chosen for their implementation, it is fundamental to judiciously choose this parameter.

Informal guidelines for the optimal scaling problem have been proposed by, among others, [3] and [4], but the first theoretical results were obtained by [11]. In particular, the authors considered $d$-dimensional target distributions with i.i.d. components and studied the asymptotic behavior (as $d \to \infty$) of Metropolis algorithms with Gaussian proposals. It was proven that under some regularity conditions for the target distribution, the asymptotic acceptance rate should be tuned to be approximately 0.234 for optimal per-

Received January 2006; revised December 2006.
[1]Supported by NSERC of Canada.

formance of the algorithm. It was also shown that the correct proposal scaling is of the form $\ell^2/d$ for some constant $\ell$ as $d \to \infty$. The simplicity of the obtained asymptotically optimal acceptance rate (AOAR) makes these theoretical results extremely useful in practice. Optimal scaling issues have been explored by other authors, namely [5, 6, 10, 12, 13].

In this paper, we carry out a similar study for $d$-dimensional target distributions with independent components. The particularity of our model is that the scaling term of each component is allowed to depend on the dimension of the target distribution, which constitutes a critical distinction with the i.i.d. case. We provide a condition under which the algorithm admits the same limiting diffusion process and the same AOAR as those found in [11]. This is achieved, in the first place, by determining the appropriate form for the proposal scaling as a function of $d$, which is now different from the i.i.d. case. Then, by verifying $\mathcal{L}^1$ convergence of generators, we prove that the sequence of stochastic processes formed by, say, the $i^*$th component of each Markov chain (appropriately rescaled) converges to a Langevin diffusion process with a certain speed measure. Obtaining the AOAR is thus a simple matter of optimizing the speed measure of the diffusion.

The paper is structured as follows. In Section 2, we describe the Metropolis algorithm and introduce the target distribution setting. The main results are presented in Section 3, along with a discussion concerning inhomogeneous proposal distributions and some extensions. We prove the theorems in Section 4 using lemmas proved in Sections 5 and 6, finally concluding the paper with a discussion.

## 2. Sampling from the target distribution.

2.1. *The Metropolis algorithm.* The idea behind the Metropolis algorithm is to generate a Markov chain $\mathbf{X}_0, \mathbf{X}_1, \ldots$ having the target distribution as a stationary distribution. In particular, suppose that $\pi$ is a $d$-dimensional probability density of interest with respect to Lebesgue measure. Also, let the proposed moves be normally distributed around $\mathbf{x}$, that is, $N(\mathbf{x}, \sigma^2 I_d)$ for some $\sigma^2$ and where $I_d$ is the $d$-dimensional identity matrix. The Metropolis algorithm thus proceeds as follows. Given $\mathbf{X}_t$, the state of the chain at time $t$, a value $\mathbf{Y}_{t+1}$ is generated from the normal density $q(\mathbf{X}_t, \mathbf{y}) \, d\mathbf{y}$. The probability of accepting the proposed value $\mathbf{Y}_{t+1}$ as the new value for the chain is $\alpha(\mathbf{X}_t, \mathbf{Y}_{t+1})$, where

$$\alpha(\mathbf{x}, \mathbf{y}) = \begin{cases} \min\left(1, \dfrac{\pi(\mathbf{y})}{\pi(\mathbf{x})}\right), & \pi(\mathbf{x})q(\mathbf{x}, \mathbf{y}) > 0, \\ 1, & \pi(\mathbf{x})q(\mathbf{x}, \mathbf{y}) = 0. \end{cases}$$

If the proposed move is accepted, the chain jumps to $\mathbf{X}_{t+1} = \mathbf{Y}_{t+1}$; otherwise, it stays where it is and $\mathbf{X}_{t+1} = \mathbf{X}_t$.

In order to have some level of optimality in the performance of the algorithm, care must be exercised when choosing $\sigma^2$. If it is too small, the proposed jumps will be too short and in spite of a very high acceptance rate, simulation will move very slowly to the target distribution. At the opposite extreme, a large scaling value will generate jumps in low target density regions, resulting in the rejection of the proposed moves and a chain that stands still most of the time.

Before finding an appropriate value for $\sigma^2$ between these extremes, we first define a criterion which is closely related to the algorithm efficiency. The notion of $\pi$-average acceptance rate is defined in [11] as $\mathrm{E}[1 \wedge \frac{\pi(\mathbf{Y})}{\pi(\mathbf{X})}] = \iint \pi(\mathbf{x})\alpha(\mathbf{x},\mathbf{y})q(\mathbf{x},\mathbf{y})\,d\mathbf{x}\,d\mathbf{y}$ for the $d$-dimensional Metropolis algorithm.

2.2. *The target distribution.* Consider the following $d$-dimensional target density

$$\pi(d, \mathbf{x}^{(d)}) = \prod_{j=1}^{d} \theta_j(d) f(\theta_j(d) x_j). \tag{1}$$

We impose the following regularity conditions on the density $f$: $f$ is a positive $C^2$ function and $(\log f(x))'$ is Lipschitz continuous. We also suppose that $\mathrm{E}[(\frac{f'(X)}{f(X)})^4] = \int_{\mathbf{R}} (\frac{f'(x)}{f(x)})^4 f(x)\,dx < \infty$ and $\mathrm{E}[(\frac{f''(X)}{f(X)})^2] < \infty$.

The $d$ target components, although independent, are, however, not identically distributed. In particular, we consider the case where the scaling terms $\theta_j^{-2}(d)$, $j=1,\ldots,d$, take the following form:

$$\mathbf{\Theta}^{-2}(d) = \bigg(\frac{K_1}{d^{\lambda_1}},\ldots,\frac{K_n}{d^{\lambda_n}},\underbrace{\frac{K_{n+1}}{d^{\gamma_1}},\ldots,\frac{K_{n+1}}{d^{\gamma_1}}}_{c(\mathcal{J}(1,d))},\ldots,\underbrace{\frac{K_{n+m}}{d^{\gamma_m}},\ldots,\frac{K_{n+m}}{d^{\gamma_m}}}_{c(\mathcal{J}(m,d))}\bigg).$$

Ultimately, we shall be interested in the limit of the target distribution as $d \to \infty$. Let $n < \infty$ denote the number of components whose scaling term appears a finite number of times in the limit of $\mathbf{\Theta}^{-2}(d)$. Also, let the $j$th of these $n$ scaling terms be $K_j/d^{\lambda_j}$, $j=1,\ldots,n$, where $\lambda_j \in (-\infty,\infty)$ and $K_j$ is some positive and finite constant. Similarly, let $0 < m < \infty$ denote the number of different scaling terms appearing infinitely often in the limit. These $m$ scaling terms are taken to be $K_{n+i}/d^{\gamma_i}$, $i=1,\ldots,m$, with $\gamma_i \in (-\infty,\infty)$. For now, we assume the constants $0 < K_{n+i} < \infty$ to be the same for all scaling terms within each of the $m$ groups. We shall relax this assumption in Section 3.2.

For $i = 1,\ldots,m$, define the sets $\mathcal{J}(i,d) = \{j \in \{1,\ldots,d\}; \theta_j^{-2}(d) = \frac{K_{n+i}}{d^{\gamma_i}}\}$. The $i$th set thus contains positions of components with a scaling term equal to $K_{n+i}/d^{\gamma_i}$. These sets are such that $\dot{\bigcup}_{i=1}^{m} \mathcal{J}(i,d) = \{n+1,\ldots,d\}$.

Since each of the $m$ groups of scaling terms might occupy different proportions of $\boldsymbol{\Theta}^{-2}(d)$, we also define the cardinality of the sets $\mathcal{J}(i,d)$:

$$（2）\quad c(\mathcal{J}(i,d)) = \#\left\{j \in \{1,\ldots,d\}; \theta_j^{-2}(d) = \frac{K_{n+i}}{d^{\gamma_i}}\right\}, \qquad i = 1,\ldots,m,$$

where $c(\mathcal{J}(i,d))$ is assumed to be some polynomial function of the dimension satisfying $\lim_{d\to\infty} c(\mathcal{J}(i,d)) = \infty$.

It will be convenient to rearrange the terms of $\boldsymbol{\Theta}^{-2}(d)$ so that all of the different scaling terms appear at one of the first $n+m$ positions:

$$（3）\quad \boldsymbol{\Theta}^{-2}(d) = \bigg(\frac{K_1}{d^{\lambda_1}},\ldots,\frac{K_n}{d^{\lambda_n}},\frac{K_{n+1}}{d^{\gamma_1}},\ldots,\frac{K_{n+m}}{d^{\gamma_m}},$$
$$\frac{K_{n+1}}{d^{\gamma_1}},\ldots,\frac{K_{n+m}}{d^{\gamma_m}},\ldots,\frac{K_{n+1}}{d^{\gamma_1}},\ldots,\frac{K_{n+m}}{d^{\gamma_m}}\bigg).$$

This helps to identify each component being studied as $d \to \infty$ without referring to a component that would otherwise be at an infinite position.

Without loss of generality, we assume the first $n$ and the next $m$ scaling terms in (3) to be respectively arranged according to an asymptotic increasing order, in the following sense. If $\preceq$ means "is asymptotically smaller than or equal to," then we have $\theta_1^{-2}(d) \preceq \cdots \preceq \theta_n^{-2}(d)$ and similarly $\theta_{n+1}^{-2}(d) \preceq \cdots \preceq \theta_{n+m}^{-2}(d)$, which respectively implies that $-\infty < \lambda_n \leq \lambda_{n-1} \leq \cdots \leq \lambda_1 < \infty$ and $-\infty < \gamma_m \leq \gamma_{m-1} \leq \cdots \leq \gamma_1 < \infty$. Based on this ordering, the asymptotically smallest scaling term obviously has to be either $\theta_1^{-2}(d)$ or $\theta_{n+1}^{-2}(d)$.

Our goal is to study the limiting distribution of each component forming the $d$-dimensional Markov process. To this end, we set the scaling term of the target component of interest equal to 1 $[\theta_{i^*}(d) = 1]$. This adjustment, necessary to obtain a nontrivial limiting process, is performed without loss of generality by applying a linear transformation to the target distribution. In particular, when the first component of the chain is studied ($i^* = 1$), we set $\theta_1^{-2}(d) = 1$ and adjust the other scaling terms accordingly. $\boldsymbol{\Theta}^{-2}(d)$ thus varies according to the component of interest $i^*$ considered.

2.3. *The proposal distribution and its scaling.* A crucial step in the implementation of Metropolis algorithms is to determine the optimal form for the proposal scaling as a function of $d$. Intuitively, it makes sense that $\sigma^2(d)$ depends on the asymptotically smallest scaling term in $\boldsymbol{\Theta}^{-2}(d)$. Otherwise, the proposed moves might be too large for the components with smaller scaling terms, resulting in a high rejection rate and compromising the convergence of the algorithm.

Moreover, as the dimension of the target increases, more individual moves are proposed in a single step and it is thus more likely to generate an improbable move for one of the components. To rectify the situation, it is recommended to decrease the proposal scaling as a function of $d$.

Hence, the optimal form of the proposal scaling turns out to be $\sigma^2(d) = \ell^2/d^\alpha$, where $\ell^2$ is some constant and $\alpha$ is the smallest number satisfying

$$(4) \qquad \lim_{d \to \infty} \frac{d^{\lambda_1}}{d^\alpha} < \infty \quad \text{and} \quad \lim_{d \to \infty} \frac{d^{\gamma_i} c(\mathcal{J}(i,d))}{d^\alpha} < \infty, \qquad i = 1, \ldots, m.$$

Therefore, at least one of these $m + 1$ limits converges to some positive constant, while the other ones converge to 0. Since the scaling term of the component studied is taken to be 1, the largest possible form for the proposal scaling is $\sigma^2 = \sigma^2(d) = \ell^2$ and so it never diverges as $d$ grows.

By its nature, the Metropolis algorithm is a discrete-time process. Since space (the proposal scaling) is a function of the dimension of the target distribution, we also have to rescale the time between each step in order to obtain a nontrivial limiting process as $d \to \infty$.

Let $\mathbf{Z}^{(d)}(t)$ be the time-$t$ value of the process sped up by a factor of $d^\alpha$; in particular, $\mathbf{Z}^{(d)}(t) = (X_1^{(d)}([d^\alpha t]), \ldots, X_d^{(d)}([d^\alpha t]))$, where $[\cdot]$ is the "integer part" function. Instead of proposing only one move, the sped-up process has the possibility of moving on average $d^\alpha$ times during each unit time interval. We are now ready to study the limiting behavior of every component of the sequence of processes $\{\mathbf{Z}^{(d)}(t), t \geq 0\}$ as $d \to \infty$.

## 3. Optimizing the sampling procedure.

3.1. *Optimal value for $\ell$.* We shall now present explicit asymptotic results allowing us to optimize $\ell^2$, the constant term of $\sigma^2(d)$. We first introduce a weak convergence result for the process $\{\mathbf{Z}^{(d)}(t), t \geq 0\}$ and most importantly in practice, we transform the achieved conclusion into a statement about efficiency as a function of acceptance rate, as was done in [11].

We denote weak convergence in the Skorokhod topology by $\Rightarrow$, standard Brownian motion at time $t$ by $B(t)$ and the standard normal cumulative distribution function (c.d.f.) by $\Phi(\cdot)$. Moreover, recall that the scaling term of the component of interest $X_{i^*}$ is taken to be one $[\theta_{i^*}(d) = 1]$, which, as explained in Section 2.2, might require a linear transformation of $\mathbf{\Theta}^{-2}(d)$.

THEOREM 1. *Consider a Metropolis algorithm with proposal distribution $\mathbf{Y}^{(d)} \sim N(\mathbf{x}^{(d)}, \frac{\ell^2}{d^\alpha} I_d)$, where $\alpha$ satisfies* (4), *applied to a target density satisfying the specified conditions on $f$ as in* (1), *with $\theta_j^{-2}(d)$, $j = 1, \ldots, d$, as in* (3) *and $\theta_{i^*}(d) = 1$. Consider the $i^*$th component of the process $\{\mathbf{Z}^{(d)}(t), t \geq 0\}$, that is, $\{Z_{i^*}^{(d)}(t), t \geq 0\} = \{X_{i^*}^{(d)}([d^\alpha t]), t \geq 0\}$ and let $\mathbf{X}^{(d)}(0)$ be distributed according to the target density $\pi$ in* (1).

*We have $\{Z_{i^*}^{(d)}(t), t \geq 0\} \Rightarrow \{Z(t), t \geq 0\}$, where $Z(0)$ is distributed according to the density $f$ and $\{Z(t), t \geq 0\}$ satisfies the Langevin stochastic* *differential equation (SDE)*

$$dZ(t) = \upsilon(\ell)^{1/2} \, dB(t) + \tfrac{1}{2}\upsilon(\ell)(\log f(Z(t)))' \, dt,$$

*if and only if*

(5) $$\lim_{d \to \infty} \frac{d^{\lambda_1}}{\sum_{j=1}^{n} d^{\lambda_j} + \sum_{i=1}^{m} c(\mathcal{J}(i,d)) \, d^{\gamma_i}} = 0.$$

*Here, $\upsilon(\ell) = 2\ell^2 \Phi(-\ell\sqrt{E_R}/2)$ and*

(6) $$E_R = \lim_{d \to \infty} \sum_{i=1}^{m} \frac{c(\mathcal{J}(i,d))}{d^{\alpha}} \frac{d^{\gamma_i}}{K_{n+i}} \mathrm{E}\left[\left(\frac{f'(X)}{f(X)}\right)^2\right],$$

*with $c(\mathcal{J}(i,d))$ as in* (2).

Intuitively, when none of the target components possesses a scaling term significantly smaller than those of the other components, the limiting process is the same as that found in [11]. Note that the numerator in condition (5) is based on $\theta_1^{-2}(d)$ only, which is not necessarily the asymptotically smallest scaling term. Technically, we should then also verify that this condition is still satisfied when $\theta_1^{-2}(d)$ is replaced by $\theta_{n+1}^{-2}(d)$; however, this is ensured by the presence of the term $c(\mathcal{J}(1,d))\theta_{n+1}^2(d)$ in the denominator.

The function $\upsilon(\ell)$ is sometimes interpreted as the speed measure of the diffusion process. As this quantity is proportional to the mixing rate of the algorithm, it suffices to maximize the function $\upsilon(\ell)$ in order to optimize the efficiency of the algorithm.

Let $a(d, \ell)$ be the $\pi$-average acceptance rate defined in Section 2.1, but where the dependence on the dimension and the proposal scaling are now made explicit. The following corollary introduces the optimal value $\hat{\ell}$ and AOAR leading to greatest efficiency of the Metropolis algorithm.

COROLLARY 2. *In the setting of Theorem* 1, *we have* $\lim_{d \to \infty} a(d, \ell) = 2\Phi(-\ell\sqrt{E_R}/2) \equiv a(\ell)$. *Furthermore, $\upsilon(\ell)$ is maximized at the unique value $\hat{\ell} = 2.38/\sqrt{E_R}$ for which $a(\hat{\ell}) = 0.234$ (to three decimal places).*

For a high-dimensional target distribution as defined in Section 2.2 and having no component converging significantly faster than the others, the value $\ell$ should be chosen such that the acceptance rate is close to 0.234 in order to optimize the efficiency of the Metropolis algorithm.

Theorem 1 may be used to determine whether or not the AOAR for sampling from any multivariate normal distribution with covariance matrix $\Sigma$ is 0.234. Since normal random variables are invariant under orthogonal transformations, we can transform $\Sigma$ into a diagonal matrix where the eigenvalues of $\Sigma$ constitute the diagonal elements. The eigenvalues can then be used to determine whether or not condition (5) is satisfied and hence to determine whether or not $2.38/\sqrt{E_R}$ is the optimal scaling for the proposal distribution. For example, consider $\Sigma$ with $\sigma_i^2 = 2$, $i = 1, \ldots, d$, and $\sigma_{ij} = 1$, $j \neq i$. The $d$ eigenvalues of $\Sigma$ are $(d, 1, \ldots, 1)$ and satisfy condition (5). For a relatively high-dimensional multivariate normal with such a correlation structure, it is thus optimal to tune the acceptance rate to 0.234. Note, however, that not all $d$ components mix at the same rate. When studying any of the last $d - 1$ components, the vector $\mathbf{\Theta}^{-2}(d) = (d, 1, \ldots, 1)$ is appropriate, so $\sigma^2(d) = \ell^2/d$ and these components thus mix in $O(d)$ iterations. When studying the first component, we need to linearly transform the scaling vector so that $\theta_1^{-2}(d) = 1$. We then use $\mathbf{\Theta}^{-2}(d) = (1, 1/d, \ldots, 1/d)$, so $\sigma^2(d) = \ell^2/d^2$ and this component mixes according to $O(d^2)$.

Now, consider the simple model where $X_1 \sim N(0, 1)$ and $X_j \sim N(X_1, 1)$ for $j = 2, \ldots, d$. The joint distribution of $\mathbf{X}^{(d)}$ is multivariate normal with mean 0 and $d \times d$ covariance matrix such that $\sigma_1^2 = 1$, $\sigma_2^2 = \cdots = \sigma_d^2 = 2$ and $\sigma_{jk} = 1$, $\forall j \neq k$. Using the $d$ eigenvalues, which are $O(d)$, $O(1/d)$ and 1 with multiplicity $d - 2$, we thus conclude that condition (5) is violated and that 0.234 might not be optimal, even though the distribution is normal [see Theorem 5 of Section 3.2 when dealing with more general $\theta_j(d)$'s].

The previous example might seem surprising as multivariate normal distributions have long been believed to behave as i.i.d. target distributions in limit. A natural question to ask, then, is what happens when condition (5) is not satisfied? In such a case, the algorithm can be shown to admit the same limiting Langevin diffusion process, but with a different speed measure. Furthermore, the AOAR is found to be smaller than the usual 0.234. For more details on this case, see [1]. For a better picture of the applicability of these results, examples and simulation studies for various statistical models are presented in [2].

3.2. *Inhomogeneous proposal scaling and extensions.* Thus far, we have assumed $\sigma^2(d) = \ell^2/d^\alpha$ to be the same for all $d$ components. It is natural to wonder whether adjusting the proposal scaling as a function of $d$ for each component would yield a better performance of the algorithm. An important point to keep in mind is that for $\{\mathbf{Z}^{(d)}(t), t \geq 0\}$ to be a stochastic process, we must speed up time by the same factor for every component. Otherwise, some components would move more frequently than others in the same time interval and since the acceptance probability of the proposed moves depends on all $d$ components, this would violate the definition of a stochastic process.

The inhomogeneous scheme we adopt is the following: we personalize the proposal scaling of the last $d - n$ components only, implying that the proposal scaling of the first $n$ components is the same as it would have been under the homogeneity assumption. We then treat each of the $m$ groups of scaling terms appearing infinitely often as a different portion of the scaling vector and determine the appropriate $\alpha$ for each group.

In particular, consider $\boldsymbol{\Theta}^{-2}(d)$ in (3) and let the proposal scaling of $X_j$ be $\sigma_j^2(d) = \ell^2/d^{\alpha_j}$, where $\alpha_j = \alpha$ for $j = 1, \ldots, n$ and $\alpha_j$ is the smallest value such that $\lim_{d\to\infty} c(\mathcal{J}(i,d))\, d^{\gamma_i}/d^{\alpha_j} < \infty$ for $j = n+1, \ldots, d$, $j \in \mathcal{J}(i,d)$. In order to study the component $X_{i^*}$, we still assume that $\theta_{i^*}(d) = 1$, but we now let $\mathbf{Z}^{(d)}(t) = \mathbf{X}^{(d)}([d^{\alpha_{i^*}} t])$. We have the following result.

THEOREM 3. *In the setting of Theorem 1, but with the proposal scaling as just described, the conclusions of Theorem 1 and Corollary 2 are preserved and $E_R$ is now expressed as*

$$E_R = \lim_{d\to\infty} \sum_{i=1}^{m} \frac{c(\mathcal{J}(i,d))}{d^{\alpha_{n+i}}} \frac{d^{\gamma_i}}{K_{n+i}} \mathrm{E}\left[\left(\frac{f'(X)}{f(X)}\right)^2\right].$$

Since the proposal scaling is now adjusted to suit every distinct group of components, each constant term $K_{n+1}, \ldots, K_{n+m}$ has an impact on the limiting process, yielding a larger value for $E_R$. Hence, the optimal value $\hat{\ell} = 2.38/\sqrt{E_R}$ is now smaller than with homogeneous proposal scaling. When the proposal scaling of all components was based on $\alpha$ in Section 3.1, the algorithm had to compensate for the fact that $\alpha$ is chosen as small as possible and thus possibly too small for certain groups of components, with a larger value for $\hat{\ell}^2$.

The conclusions of Section 3.1 also extend to more general target distribution settings. First, we can relax the assumption of equality among the constant terms of $\theta_j^{-2}(d)$ for $j \in \mathcal{J}(i,d)$. In particular, let

$$\boldsymbol{\Theta}^{-2}(d) = \left(\frac{K_1}{d^{\lambda_1}}, \ldots, \frac{K_n}{d^{\lambda_n}}, \frac{K_{n+1}}{d^{\gamma_1}}, \ldots, \frac{K_{n+c(\mathcal{J}(1,d))}}{d^{\gamma_1}}, \ldots, \right.$$
(7)
$$\left. \frac{K_{n+\sum_{i=1}^{m-1} c(\mathcal{J}(i,d))+1}}{d^{\gamma_m}}, \ldots, \frac{K_d}{d^{\gamma_m}}\right).$$

We assume that $\{K_j, j \in \mathcal{J}(i,d)\}$ are i.i.d. and chosen randomly from some distribution with $\mathrm{E}[K_j^{-2}] < \infty$. Without loss of generality, we denote $\mathrm{E}[K_j^{-1}] = b_i$ for $j \in \mathcal{J}(i,d)$. Recall that the scaling term of the component of interest cannot depend on $d$, so we have $\theta_{i^*}^{-2}(d) = K_{i^*}$.

To support the previous modifications, we now suppose that $-\infty < \gamma_m < \gamma_{m-1} < \cdots < \gamma_1 < \infty$. In addition, we assume that there does not exist a $\lambda_j$, $j = 1, \ldots, n$, equal to one of the $\gamma_i$, $i = 1, \ldots, m$. This means that if there is an infinite number of scaling terms of the same order, they must necessarily belong to the same of the $m$ groups. We obtain the following result.

THEOREM 4. *Consider the setting of Theorem 1, except with $\mathbf{\Theta}^{-2}(d)$ as in (7) and $\theta_{i^*} = K_{i^*}^{-1/2}$. We have $\{Z_{i^*}^{(d)}(t), t \geq 0\} \Rightarrow \{Z(t), t \geq 0\}$, where $Z(0)$ is distributed according to the density $\theta_{i^*} f(\theta_{i^*} x)$ and $\{Z(t), t \geq 0\}$ satisfies the Langevin SDE*

$$dZ(t) = (\upsilon(\ell))^{1/2} \, dB(t) + \tfrac{1}{2}\upsilon(\ell)(\log f(\theta_{i^*} Z(t)))' \, dt,$$

*if and only if condition (5) is satisfied. Here, $\upsilon(\ell)$ is as in Theorem 1 and*

$$E_R = \lim_{d \to \infty} \sum_{i=1}^{m} \frac{c(\mathcal{J}(i,d)) \, d^{\gamma_i}}{d^{\alpha}} b_i \mathrm{E}\left[\left(\frac{f'(X)}{f(X)}\right)^2\right],$$

*with $c(\mathcal{J}(i,d)) = \#\{j \in \{n+1,\ldots,d\}; \theta_j(d) \text{ is } O(d^{\gamma_i/2})\}$. Furthermore, the conclusions of Corollary 2 are preserved.*

The previous results can also be extended to more general functions $c(\mathcal{J}(i,d))$, $i = 1,\ldots,m$, and $\theta_j(d)$, $j = 1,\ldots,d$. In order to have sensible limiting theory, however, we restrict our attention to functions for which the limit exists as $d \to \infty$. As before, we must have $c(\mathcal{J}(i,d)) \to \infty$ as $d \to \infty$. We even allow $\{\theta_j^{-2}(d), j \in \mathcal{J}(i,d)\}$ to vary within each of the $m$ groups, provided they are of the same order. That is, for $j \in \mathcal{J}(i,d)$, we suppose that $\lim_{d \to \infty} \theta_j(d)/\theta_i'(d) = K_j^{-1/2}$ for some reference function $\theta_i'(d)$ and some constant $K_j$ coming from the distribution described in Theorem 4.

As for Theorem 4, we assume that if there are infinitely many scaling terms of a certain order, then they must all belong to one of the $m$ groups. Hence, $\mathbf{\Theta}^{-2}(d)$ contains at least $m$ and at most $n+m$ functions of different orders. The positions of the elements belonging to the $i$th group are thus

$$(8) \quad \mathcal{J}(i,d) = \left\{j \in \{1,\ldots,d\}; 0 < \lim_{d \to \infty} \frac{\theta_i'^2(d)}{\theta_j^2(d)} < \infty\right\}, \qquad i \in \{1,\ldots,m\}.$$

For such target distributions, we define the proposal scaling to be $\sigma^2(d) = \ell^2 \sigma_\alpha^2(d)$, with $\sigma_\alpha^2(d)$ the function of largest possible order such that

$$(9) \quad \begin{aligned} &\lim_{d \to \infty} \theta_1^2(d)\sigma_\alpha^2(d) < \infty \qquad \text{and} \\ &\lim_{d \to \infty} c(\mathcal{J}(i,d))\theta_i'^2(d)\sigma_\alpha^2(d) < \infty, \qquad i = 1,\ldots,m. \end{aligned}$$

THEOREM 5. *Under the setting of Theorem 4, but with proposal scaling $\sigma^2(d) = \ell^2 \sigma_\alpha^2(d)$, where $\sigma_\alpha^2(d)$ satisfies (9) and with general functions for $c(\mathcal{J}(i,d))$ and $\theta_j(d)$ as defined previously, the conclusions of Theorem 4 are preserved, provided that*

$$\lim_{d \to \infty} \frac{\theta_1^2(d)}{\sum_{j=1}^{n} \theta_j^2(d) + \sum_{i=1}^{m} c(\mathcal{J}(i,d))\theta_i'^2(d)} = 0$$

*holds instead of condition* (5) *and with*

$$E_R = \lim_{d\to\infty} \sum_{i=1}^{m} c(\mathcal{J}(i,d))\theta_i'^2(d)\sigma_\alpha^2(d)b_i \mathrm{E}\bigg[\bigg(\frac{f'(X)}{f(X)}\bigg)^2\bigg],$$

*where $c(\mathcal{J}(i,d))$ is the cardinality function of* (8).

This theorem assumes quite a general form for the scaling terms of the target distribution and allows for a lot of flexibility.

**4. Proofs of theorems.** We now present the proof of Theorem 1; those of the theorems in Section 3.2 being similar, we just outline the main differences. The proofs are based on Theorem 8.2 of Chapter 4 in [7] which roughly says that for the finite-dimensional distributions of a sequence of processes to converge weakly to those of some Markov process, it is sufficient to verify $\mathcal{L}^1$ convergence of their generators. Then, Corollary 8.6 of the same chapter provides further conditions for our sequence of processes to be relatively compact and thus to reach weak convergence of the stochastic processes themselves. Specifically, it is easily verified that $C_c^\infty$, the space of infinitely differentiable functions with compact support, is an algebra that strongly separates points. Since the algorithm starts in stationarity, $\mathbf{X}^{(d)}(t) \sim \pi \,\forall t > 0$. Using a method similar to the proof of Lemma 7, we show that $\mathrm{E}[(Gh(d,\mathbf{X}^{(d)}))^2]$ is bounded by some constant for all $d \geq 1$, where $G$ is the generator of the sped-up Metropolis algorithm appearing in Section 4.2; this ensures relative compactness.

Our task is then to focus on the $\mathcal{L}^1$ convergence of the generators. To this end, we base our approach on the proof for the Metropolis algorithm case in [10]. Note, however, that the authors instead prove uniform convergence of generators, and this could not be used in the present situation.

The generator is written in terms of an arbitrary test function $h$ which can usually be any smooth function; in our case, we restrict our attention to functions in $C_c^\infty$. Since the limiting process obtained is a diffusion, it follows that $C_c^\infty$ is a core for the generator by Theorem 2.1 of Chapter 8 in [7], so instead of verifying $\mathcal{L}^1$ convergence of the generators for all functions $h$ in the domain of $G_L$, we shall be allowed to work with functions belonging to this core only.

In order to ease notation, we adopt the following convention for defining vectors: $\mathbf{X}^{(b-a)} = (X_{a+1}, \ldots, X_b)$ The minus sign appearing outside the brackets (e.g., $\mathbf{X}^{(b-a)-}$) means that the component of interest, $X_{i^*}$, is excluded. We also use the following notation for conditional expectations: $\mathrm{E}[f(X,Y)|X] = \mathrm{E}_Y[f(X,Y)]$. When there is no subscript, the expectation is taken with respect to all random variables included in the expression.

4.1. *Restrictions on the proposal scaling.* We first transform condition (5) into a statement about the proposal scaling and its parameter $\alpha$. For this condition to be satisfied, we must equivalently have

$$\lim_{d\to\infty} \frac{K_1}{d^{\lambda_1}}\left(\frac{d^{\lambda_1}}{K_1} + \cdots + \frac{d^{\lambda_n}}{K_n}\right)$$
$$+ \lim_{d\to\infty} \frac{K_1}{d^{\lambda_1}}\left(c(\mathcal{J}(1,d))\frac{d^{\gamma_1}}{K_{n+1}} + \cdots + c(\mathcal{J}(m,d))\frac{d^{\gamma_m}}{K_{n+m}}\right) = \infty.$$

Since the first term on the left-hand side is finite, there is at least one $i \in \{1,\ldots,m\}$ such that $\lim_{d\to\infty} \theta_1^{-2}(d)c(\mathcal{J}(i,d))\frac{d^{\gamma_i}}{K_{n+i}} = \infty$. Consequently, the choice of $\alpha$ in (4) must be based on one of the groups of scaling terms appearing infinitely often. If we had $\alpha = \lambda_1$, this would mean that $\lim_{d\to\infty} \frac{c(\mathcal{J}(i,d))\,d^{\gamma_i}}{d^\alpha} = \infty$ for all $i$ for which the previous limit was diverging, which contradicts the definition of $\alpha$. When condition (5) is satisfied, it thus follows that $\lim_{d\to\infty} d^{\lambda_1}/d^\alpha = 0$ and $\theta_1^{-2}(d)$ does not govern $\alpha$; the parameter $\alpha$ is then strictly greater than 0, regardless of which component is under consideration.

4.2. *Proof of Theorem* 1. For an arbitrary test function $h \in C_c^\infty$, we show that

$$\lim_{d\to\infty} \mathrm{E}[|Gh(d,\mathbf{X}^{(d)}) - G_L h(X_{i^*})|] = 0,$$

where $Gh(d,\mathbf{X}^{(d)}) = d^\alpha \mathrm{E}_{\mathbf{Y}^{(d)}}[(h(Y_{i^*}) - h(X_{i^*}))(1 \wedge (\pi(d,\mathbf{Y}^{(d)})/\pi(d,\mathbf{X}^{(d)})))]$ is the discrete-time generator of the sped-up Metropolis algorithm and $G_L h(X_{i^*}) = v(\ell)[\frac{1}{2}h''(X_{i^*}) + \frac{1}{2}h'(X_{i^*})(\log f(X_{i^*}))']$ is the generator of a Langevin diffusion process with speed measure $v(\ell)$, as in Theorem 1.

According to Lemma 7, we have $\lim_{d\to\infty} \mathrm{E}[|Gh(d,\mathbf{X}^{(d)}) - \widetilde{G}h(d,\mathbf{X}^{(d)})|] = 0$, where

$$\widetilde{G}h(d,\mathbf{X}^{(d)})$$
$$= \tfrac{1}{2}\ell^2 h''(X_{i^*})\mathrm{E}_{\mathbf{Y}^{(d)-}}\left[1 \wedge \exp\left\{\sum_{j=1,j\neq i^*}^d \varepsilon(d,X_j,Y_j)\right\}\right]$$
$$+ \ell^2 h'(X_{i^*})(\log f(X_{i^*}))'\mathrm{E}_{\mathbf{Y}^{(d)-}}$$
$$\times \left[\exp\left\{\sum_{j=1,j\neq i^*}^d \varepsilon(d,X_j,Y_j)\right\}; \sum_{j=1,j\neq i^*}^d \varepsilon(d,X_j,Y_j) < 0\right]$$

and $\varepsilon(d,X_j,Y_j)$ is as in (10). To prove the theorem, we are thus left to show $\mathcal{L}^1$ convergence of the generator $\widetilde{G}h(d,\mathbf{X}^{(d)})$ to the generator of the Langevin diffusion.

Substituting explicit expressions for the generators, grouping some terms and using the triangle inequality, we obtain

$$\mathrm{E}[|\widetilde{G}h(d, \mathbf{X}^{(d)}) - G_L h(X_{i^*})|]$$
$$\leq \ell^2 \mathrm{E}_{\mathbf{X}^{(d)-}} \left[ \left| \frac{1}{2} \mathrm{E}_{\mathbf{Y}^{(d)-}} \left[ 1 \wedge \exp\left\{ \sum_{j=1, j \neq i^*}^{d} \varepsilon(d, X_j, Y_j) \right\} \right] \right. \right.$$
$$\left. \left. - \Phi\left( -\frac{\ell \sqrt{E_R}}{2} \right) \right| \right] \mathrm{E}[|h''(X_{i^*})|]$$
$$+ \ell^2 \mathrm{E}_{\mathbf{X}^{(d)-}} \left[ \left| \mathrm{E}_{\mathbf{Y}^{(d)-}} \left[ \exp\left\{ \sum_{j=1, j \neq i^*}^{d} \varepsilon(d, X_j, Y_j) \right\}; \right. \right. \right.$$
$$\left. \left. \sum_{j=1, j \neq i^*}^{d} \varepsilon(d, X_j, Y_j) < 0 \right] - \Phi\left( -\frac{\ell \sqrt{E_R}}{2} \right) \right| \right]$$
$$\times \mathrm{E}[|h'(X_{i^*})(\log f(X_{i^*}))'|].$$

Since the function $h$ has compact support, it follows that $h$ itself and its derivatives are bounded in absolute value by some constant. As a result, $\mathrm{E}[|h''(X_{i^*})|]$ and $\mathrm{E}[|h'(X_{i^*})(\log f(X_{i^*}))'|]$ are both bounded by $K$, say. Using Lemmas 8 and 9, we then conclude that the first expectation on the right-hand side goes to 0 as $d \to \infty$; we reach the same conclusion for the first expectation of the second term by applying Lemmas 10 and 11.

4.3. *Proof of Theorem* 4. The main difference with the proof of Theorem 1 occurs when working with the $m$ groups formed of infinitely many components. Since the constant terms are now random, we cannot factorize the scaling terms of components belonging to the same group. However, this difficulty is easily overcome by changes of variable and the use of conditional expectations; for instance, a typical quantity we must work with is

$$\frac{1}{d^\alpha} \sum_{j \in \mathcal{J}(i,d)} \left( \frac{d}{dX_j} \log \theta_j(d) f(\theta_j(d) X_j) \right)^2$$
$$= \frac{c(\mathcal{J}(i,d)) d^{\gamma_i}}{d^\alpha} \left[ \frac{1}{c(\mathcal{J}(i,d))} \sum_{j \in \mathcal{J}(i,d)} \left( \frac{d}{dX_j} \log \frac{1}{\sqrt{K_j}} f\left( \frac{X_j}{\sqrt{K_j}} \right) \right)^2 \right].$$

By the weak law of large numbers (WLLN; see, e.g., [15]), the term in brackets converges to $b_i \mathrm{E}[(f'(X)/f(X))^2]$. Instead of carrying the term $\theta_{n+i}^2(d) = d^{\gamma_i}/K_{n+i}$ as before, we thus carry $b_i d^{\gamma_i}$.

4.4. *Proof of Theorem 5.* The general forms of the functions $c(\mathcal{J}(i,d))$, $i=1,\ldots,m$, and $\theta_j(d)$, $j=1,\ldots,d$ necessitate a more elaborate notation, but do not affect the body of the proof. Instead, what alters the demonstration is the fact that $\theta_j(d)$ for $j \in \mathcal{J}(i,d)$ are allowed to be different functions of $d$ provided they are of the same order. Because of this particularity, we must write $\theta_j(d) = K_j^{-1/2}\theta_i'(d)\theta_j^*(d)/\theta_i'(d)$, where $\theta_j^*(d)$ is implicitly defined. We can then continue with the proof as usual, factoring the term $b_i\theta_i'(d)$ instead of $\theta_{n+i}^2(d)$ in Theorem 1 (or $b_i d^{\gamma_i}$ in Theorem 4). Since $\lim_{d\to\infty} \theta_j^*(d)/\theta_i'(d) = 1$, the rest of the proof can be repeated with minor modifications.

## 5. Equivalent generator and other results.

5.1. *Convergence of an approximation term.*

LEMMA 6. *For $i=1,\ldots,m$, let*

$$W_i(d, \mathbf{X}_{\mathcal{J}(i,d)}^{(d)}, \mathbf{Y}_{\mathcal{J}(i,d)}^{(d)}) = \frac{1}{2} \sum_{j \in \mathcal{J}(i,d)} \left(\frac{d^2}{dX_j^2} \log f(\theta_j(d)X_j)\right)(Y_j - X_j)^2$$

$$+ \frac{\ell^2}{2d^\alpha} \sum_{j \in \mathcal{J}(i,d)} \left(\frac{d}{dX_j} \log f(\theta_j(d)X_j)\right)^2,$$

*where $Y_j|X_j \sim N(X_j, \ell^2/d^\alpha)$ and $X_j$ is distributed according to the density $\theta_j(d)f(\theta_j(d)x_j)$, independently for all $j=1,\ldots,d$. Then, for $i=1,\ldots,m$,*

$$\mathrm{E}_{\mathbf{Y}_{\mathcal{J}(i,d)}^{(d)}}[|W_i(d, \mathbf{X}_{\mathcal{J}(i,d)}^{(d)}, \mathbf{Y}_{\mathcal{J}(i,d)}^{(d)})|] \xrightarrow{p} 0 \qquad \text{as } d \to \infty.$$

PROOF. By Jensen's inequality, $\mathrm{E}[|W|] \leq \sqrt{\mathrm{E}[W^2]}$. Developing the square and taking the expectation conditional on $\mathbf{X}_{\mathcal{J}(i,d)}^{(d)}$, we obtain

$$\mathrm{E}_{\mathbf{Y}_{\mathcal{J}(i,d)}^{(d)}}[W_i^2(d, \mathbf{X}_{\mathcal{J}(i,d)}^{(d)}, \mathbf{Y}_{\mathcal{J}(i,d)}^{(d)})]$$

$$= \frac{\ell^4}{2d^{2\alpha}} \sum_{j \in \mathcal{J}(i,d)} \left(\frac{d^2}{dX_j^2} \log f(\theta_j(d)X_j)\right)^2$$

$$+ \frac{\ell^4}{4d^{2\alpha}} \left\{\sum_{j \in \mathcal{J}(i,d)} \left(\frac{d^2}{dX_j^2} \log f(\theta_j(d)X_j) + \left(\frac{d}{dX_j} \log f(\theta_j(d)X_j)\right)^2\right)\right\}^2.$$

Using changes of variable, we obtain

$$\mathrm{E}_{\mathbf{Y}_{\mathcal{J}(i,d)}^{(d)}}[|W_i(d, \mathbf{X}_{\mathcal{J}(i,d)}^{(d)}, \mathbf{Y}_{\mathcal{J}(i,d)}^{(d)})|]$$

$$\leq \frac{\ell^2}{\sqrt{2}d^\alpha}\theta_{n+i}^2(d)\sqrt{c(\mathcal{J}(i,d))}$$

$$\times \left(\frac{1}{c(\mathcal{J}(i,d))}\sum_{j\in\mathcal{J}(i,d)}\left(\frac{d^2}{dX_j^2}\log f(X_j)\right)^2\right)^{1/2}$$

$$+\frac{\ell^2}{2d^\alpha}\theta_{n+i}^2(d)c(\mathcal{J}(i,d))$$

$$\times\left|\frac{1}{c(\mathcal{J}(i,d))}\sum_{j\in\mathcal{J}(i,d)}\left(\frac{d^2}{dX_j^2}\log f(X_j)+\left(\frac{d}{dX_j}\log f(X_j)\right)^2\right)\right|.$$

By the WLLN, the term in parentheses on the second line converges in probability to $\mathrm{E}[(\frac{d^2}{dX^2}\log f(X))^2]$ as $d\to\infty$. Since $d^\alpha > d^{\gamma_i}\sqrt{c(\mathcal{J}(i,d))}$ and the previous expectation is bounded by some constant, the first term converges to 0 as $d\to\infty$. Given that $\theta_{n+i}^2(d)c(\mathcal{J}(i,d))/d^\alpha$ is $O(1)$ for at least one $i\in\{1,\ldots,m\}$, we must also show that the term between absolute values converges to 0. From Lemma A.1, we know that $f'(x)\to 0$ as $x\to\pm\infty$; hence, we have $\mathrm{E}[\frac{d^2}{dX_j^2}\log f(X_j)+(\frac{d}{dX_j}\log f(X_j))^2]=\int f''(x)\,dx=0$ and as $d\to\infty$, we conclude (by the WLLN) that

$$\left|\frac{1}{c(\mathcal{J}(i,d))}\sum_{j\in\mathcal{J}(i,d)}\left(\frac{d^2}{dX_j^2}\log f(X_j)+\left(\frac{d}{dX_j}\log f(X_j)\right)^2\right)\right|\xrightarrow{p}0. \qquad \square$$

5.2. *Convergence to the equivalent generator $\widetilde{G}h(d,\mathbf{X}^{(d)})$.*

LEMMA 7. *For any function $h\in C_c^\infty$, let*

$$\widetilde{G}h(d,\mathbf{X}^{(d)})$$

$$=\tfrac{1}{2}\ell^2 h''(X_{i^*})\mathrm{E}_{\mathbf{Y}^{(d)-}}\left[1\wedge\exp\left\{\sum_{j=1,j\neq i^*}^d \varepsilon(d,X_j,Y_j)\right\}\right]$$

$$+\ell^2 h'(X_{i^*})(\log f(X_{i^*}))'$$

$$\times\mathrm{E}_{\mathbf{Y}^{(d)-}}\left[\exp\left\{\sum_{j=1,j\neq i^*}^d\varepsilon(d,X_j,Y_j)\right\};\sum_{j=1,j\neq i^*}^d\varepsilon(d,X_j,Y_j)<0\right],$$

*where*

(10) $$\varepsilon(d,X_j,Y_j)=\log\frac{f(\theta_j(d)Y_j)}{f(\theta_j(d)X_j)}.$$

*If $\alpha>0$ is as defined in (4), then $\lim_{d\to\infty}\mathrm{E}[|Gh(d,\mathbf{X}^{(d)})-\widetilde{G}h(d,\mathbf{X}^{(d)})|]=0$.*

PROOF. The proof being similar to that of Lemmas A.2 and A.3 in [10], we shall omit some details. The generator of the sped-up Metropolis algorithm can be expressed as

$$Gh(d, \mathbf{X}^{(d)}) = d^\alpha \mathrm{E}_{Y_{i^*}}\left[(h(Y_{i^*}) - h(X_{i^*}))\mathrm{E}_{\mathbf{Y}^{(d)-}}\left[1 \wedge \exp\left\{\sum_{j=1}^d \varepsilon(d, X_j, Y_j)\right\}\right]\right]. \tag{11}$$

We can reexpress the inner expectation using a Taylor expansion of the minimum function with respect to $Y_{i^*}$ around $X_{i^*}$. As mentioned in [10], the generator becomes

$$Gh(d, \mathbf{X}^{(d)})$$
$$= d^\alpha \mathrm{E}_{Y_{i^*}}[(h(Y_{i^*}) - h(X_{i^*}))]\mathrm{E}_{\mathbf{Y}^{(d)-}}\left[1 \wedge \exp\left\{\sum_{j=1, j\neq i^*}^d \varepsilon(d, X_j, Y_j)\right\}\right]$$
$$+ d^\alpha (\log f(X_{i^*}))' \mathrm{E}_{Y_{i^*}}[(h(Y_{i^*}) - h(X_{i^*}))(Y_{i^*} - X_{i^*})]$$
$$\times \mathrm{E}_{\mathbf{Y}^{(d)-}}\left[\exp\left\{\sum_{j=1, j\neq i^*}^d \varepsilon(d, X_j, Y_j)\right\}; \sum_{j=1, j\neq i^*}^d \varepsilon(d, X_j, Y_j) < 0\right]$$
$$+ \frac{d^\alpha}{2}\mathrm{E}_{Y_{i^*}}[(h(Y_{i^*}) - h(X_{i^*}))(Y_{i^*} - X_{i^*})^2((\log f(U_{i^*}))')^2$$
$$\times \mathrm{E}_{\mathbf{Y}^{(d)-}}[e^{g(U_{i^*})}; g(U_{i^*}) < 0]]$$
$$+ \frac{d^\alpha}{2}\mathrm{E}_{Y_{i^*}}[(h(Y_{i^*}) - h(X_{i^*}))(Y_{i^*} - X_{i^*})^2(\log f(U_{i^*}))''$$
$$\times \mathrm{E}_{\mathbf{Y}^{(d)-}}[e^{g(U_{i^*})}; g(U_{i^*}) < 0]],$$

where $g(U_{i^*}) = \varepsilon(X_{i^*}, U_{i^*}) + \sum_{j=1, j\neq i^*}^d \varepsilon(d, X_j, Y_j)$ for some $U_{i^*} \in (X_{i^*}, Y_{i^*})$ or $(Y_{i^*}, X_{i^*})$.

We first note that all expectations computed with respect to $\mathbf{Y}^{(d)-}$ are bounded by 1, $|(\log f(U_{i^*}))''|$ is bounded by a constant and $|(\log f(U_{i^*}))'| \leq |(\log f(X_{i^*}))'| + K|Y_{i^*} - X_{i^*}|$ for some $K > 0$. Expressing $h(Y_{i^*}) - h(X_{i^*})$ as a three-term Taylor expansion and using the fact that $h$ has compact support, we can bound the expectations taken with respect to $Y_{i^*}$ and obtain

$$|Gh(d, \mathbf{X}^{(d)}) - \widetilde{G}h(d, \mathbf{X}^{(d)})|$$
$$\leq K\left(\frac{\ell^3}{d^{\alpha/2}} + \frac{\ell^4}{d^\alpha} + \frac{\ell^5}{d^{3\alpha/2}}\right)((\log f(X_{i^*}))')^2$$
$$+ K\left(\frac{\ell^4}{d^\alpha} + \frac{\ell^5}{d^{3\alpha/2}} + \frac{\ell^6}{d^{2\alpha}}\right)(1 + |(\log f(X_{i^*}))'|) + K\frac{\ell^3}{d^{\alpha/2}} + K\frac{\ell^7}{d^{5\alpha/2}}$$

for some constant $K > 0$. By assumption, $\mathrm{E}[((\log f(X_{i^*}))')^2] < \infty$, so it follows that $\mathrm{E}[|Gh(d, \mathbf{X}^{(d)}) - \widetilde{G}h(d, \mathbf{X}^{(d)})|]$ converges to 0 as $d \to \infty$. □

## 6. Volatility and drift of the diffusion.

6.1. *Convergence to an equivalent volatility.*

LEMMA 8. *We have*

$$\lim_{d\to\infty} \mathrm{E}_{\mathbf{X}^{(d)-}}\left[\left|\mathrm{E}_{\mathbf{Y}^{(d)-}}\left[1\wedge\exp\left\{\sum_{j=1,j\neq i^*}^{d}\varepsilon(d,X_j,Y_j)\right\}\right]\right.\right.$$
$$\left.\left.-\mathrm{E}_{\mathbf{Y}^{(d)-}}\left[1\wedge\exp\left\{z(d,\mathbf{Y}^{(d)-},\mathbf{X}^{(d)-})\right\}\right]\right|\right]=0,$$

*where* $\varepsilon(d,X_j,Y_j)$ *is as in* (10) *and*

$$z(d,\mathbf{Y}^{(d)-},\mathbf{X}^{(d)-})$$
$$= \sum_{j=1,j\neq i^*}^{n} \varepsilon(d,X_j,Y_j)$$

(12)
$$+ \sum_{i=1}^{m} \sum_{j\in\mathcal{J}(i,d),j\neq i^*} \frac{d}{dX_j}\log f(\theta_j(d)X_j)(Y_j-X_j)$$

$$- \frac{\ell^2}{2d^\alpha} \sum_{i=1}^{m} \sum_{j\in\mathcal{J}(i,d),j\neq i^*} \left(\frac{d}{dX_j}\log f(\theta_j(d)X_j)\right)^2.$$

PROOF. Using a Taylor expansion with three terms, we obtain

$$\mathrm{E}_{\mathbf{Y}^{(d)-}}\left[1\wedge\exp\left\{\sum_{j=1,j\neq i^*}^{d}\varepsilon(d,X_j,Y_j)\right\}\right]$$

$$= \mathrm{E}_{\mathbf{Y}^{(d)-}}\left[1\wedge\exp\left\{\sum_{j=1,j\neq i^*}^{n}\varepsilon(d,X_j,Y_j)\right.\right.$$

(13)
$$+ \sum_{i=1}^{m} \sum_{\substack{j\in\mathcal{J}(i,d)\\j\neq i^*}} \left[\frac{d}{dX_j}\log f(\theta_j(d)X_j)(Y_j-X_j)\right.$$

$$+ \frac{1}{2}\frac{d^2}{dX_j^2}\log f(\theta_j(d)X_j)(Y_j-X_j)^2$$

$$\left.\left.\left.+ \frac{1}{6}\frac{d^3}{dU_j^3}\log f(\theta_j(d)U_j)(Y_j-X_j)^3\right]\right\}\right]$$

for some $U_j \in (X_j,Y_j)$ or $(Y_j,X_j)$.

By the triangle inequality, the Lipschitz property of the function $1 \wedge e^x$ (see Proposition 2.2 in [11]) and the observation that the first two terms of the function $z(d, \mathbf{Y}^{(d)-}, \mathbf{X}^{(d)-})$ cancel the first two terms of the exponential function in (13), we get

$$\left| \mathrm{E}_{\mathbf{Y}^{(d)-}}\left[ 1 \wedge \exp\left\{ \sum_{j=1, j\neq i^*}^{d} \varepsilon(d, X_j, Y_j) \right\} \right] \right.$$
$$\left. - \mathrm{E}_{\mathbf{Y}^{(d)-}}[1 \wedge \exp\{z(d, \mathbf{Y}^{(d-)}, \mathbf{X}^{(d-)})\}] \right|$$
$$\leq \sum_{i=1}^{m} \mathrm{E}_{\mathbf{Y}^{(d)-}_{\mathcal{J}(i,d)}}[|W_i(d, \mathbf{X}^{(d)-}_{\mathcal{J}(i,d)}, \mathbf{Y}^{(d)-}_{\mathcal{J}(i,d)})|] + \sum_{i=1}^{m} c(\mathcal{J}(i,d))\ell^3 K \frac{d^{3\gamma_i/2}}{d^{3\alpha/2}}.$$

By Lemma 6, the right-hand side converges in probability to 0 as $d \to \infty$. We then apply the bounded convergence theorem to complete the proof of the lemma. $\square$

6.2. *Simplified expression for the equivalent volatility.*

LEMMA 9. *If condition* (5) *is satisfied, then*

$$\lim_{d\to\infty} \mathrm{E}_{\mathbf{X}^{(d)-}}\left[ \left| \mathrm{E}_{\mathbf{Y}^{(d)-}}[1 \wedge \exp\{z(d, \mathbf{Y}^{(d)-}, \mathbf{X}^{(d)-})\}] - 2\Phi\left(-\frac{\ell\sqrt{E_R}}{2}\right) \right| \right] = 0,$$

*where* $z(d, \mathbf{Y}^{(d)-}, \mathbf{X}^{(d)-})$ *and* $E_R$ *are as in* (12) *and* (6), *respectively.*

PROOF. For each group of components whose scaling term appears infinitely often in the limit, that is, for $i = 1, \ldots, m$, let

$$(14) \qquad R_i(d, \mathbf{x}^{(d)-}_{\mathcal{J}(i,d)}) = \frac{1}{d^\alpha} \sum_{j \in \mathcal{J}(i,d), j\neq i^*} \left( \frac{d}{dx_j} \log f(\theta_j(d) x_j) \right)^2.$$

Since $(Y_j - X_j)|X_j \sim$ i.i.d. $N(0, \ell^2/d^\alpha)$ for $j = 1, \ldots, d$, it follows that

$$z(d, \mathbf{Y}^{(d)-}, \mathbf{X}^{(d)-})|\mathbf{Y}^{(n)-}, \mathbf{X}^{(d)-}$$
$$\sim N\left( \sum_{j=1, j\neq i^*}^{n} \varepsilon(d, X_j, Y_j) - \frac{\ell^2}{2} \sum_{i=1}^{m} R_i(d, \mathbf{X}^{(d)-}_{\mathcal{J}(i,d)}), \ell^2 \sum_{i=1}^{m} R_i(d, \mathbf{X}^{(d)-}_{\mathcal{J}(i,d)}) \right).$$

Applying Proposition 2.4 in [11] allows us to obtain an expression in terms of $\Phi(\cdot)$, the c.d.f. of a standard normal random variable,

$$\mathrm{E}_{\mathbf{Y}^{(d)-}}[1 \wedge \exp\{z(d, \mathbf{Y}^{(d)-}, \mathbf{X}^{(d)-})\}]$$
$$= \mathrm{E}_{\mathbf{Y}^{(n)-}}\left[ \exp\left( \sum_{j=1, j\neq i^*}^{n} \varepsilon(d, X_j, Y_j) \right) \right.$$

$$\times \Phi\left(\frac{-\sum_{j=1,j\neq i^*}^n \varepsilon(d,X_j,Y_j) - \frac{\ell^2}{2}\sum_{i=1}^m R_i(d,\mathbf{X}_{\mathcal{J}(i,d)}^{(d)-})}{\sqrt{\ell^2 \sum_{i=1}^m R_i(d,\mathbf{X}_{\mathcal{J}(i,d)}^{(d)-})}}\right)$$

$$+ \Phi\left(\frac{\sum_{j=1,j\neq i^*}^n \varepsilon(d,X_j,Y_j) - \frac{\ell^2}{2}\sum_{i=1}^m R_i(d,\mathbf{X}_{\mathcal{J}(i,d)}^{(d)-})}{\sqrt{\ell^2 \sum_{i=1}^m R_i(d,\mathbf{X}_{\mathcal{J}(i,d)}^{(d)-})}}\right)\Bigg].$$

We note that $E_R > 0$ since there is at least one $i \in \{1,\ldots,m\}$ such that $\lim_{d\to\infty} c(\mathcal{J}(i,d))\, d^{\gamma_i}/d^\alpha > 0$. Using Propositions A.2 and A.3 and then applying Slutsky's theorem and the continuous mapping theorem, we conclude that $\exp(\sum_{j=1,j\neq i^*}^n \varepsilon(d,X_j,Y_j)) \xrightarrow{p} 1$ and

$$\Phi\left(\frac{\pm\sum_{j=1,j\neq i^*}^n \varepsilon(d,X_j,Y_j) - (\ell^2/2)\sum_{i=1}^m R_i(d,\mathbf{X}_{\mathcal{J}(i,d)}^{(d)})}{\sqrt{\ell^2 \sum_{i=1}^m R_i(d,\mathbf{X}_{\mathcal{J}(i,d)}^{(d)})}}\right)$$
$$\xrightarrow{p} \Phi\left(-\frac{\ell\sqrt{E_R}}{2}\right).$$

Since $\mathrm{E}_{\mathbf{Y}^{(d-n)-}}[1 \wedge e^{z(d,\mathbf{Y}^{(d)-},\mathbf{X}^{(d)-})}]$ is positive and bounded by 1, we use the bounded convergence theorem to conclude that $\mathrm{E}[1 \wedge e^{z(d,\mathbf{Y}^{(d)-},\mathbf{X}^{(d)-})}] \xrightarrow{p} 2\Phi(-\ell\sqrt{E_R}/2)$; we complete the proof of the lemma by reapplying the bounded convergence theorem. $\square$

### 6.3. Convergence to an equivalent drift.

LEMMA 10. *We have*

(15)
$$\lim_{d\to\infty} \mathrm{E}_{\mathbf{X}^{(d)-}}\Bigg[\Bigg| \mathrm{E}_{\mathbf{Y}^{(d)-}}\bigg[\exp\bigg\{\sum_{j=1,j\neq i^*}^d \varepsilon(d,X_j,Y_j)\bigg\};$$
$$\sum_{j=1,j\neq i^*}^d \varepsilon(d,X_j,Y_j) < 0\bigg]$$
$$- \mathrm{E}_{\mathbf{Y}^{(d)-}}[\exp\{z(d,\mathbf{Y}^{(d)-},\mathbf{X}^{(d)-})\};$$
$$z(d,\mathbf{Y}^{(d)-},\mathbf{X}^{(d)-}) < 0]\Bigg|\Bigg] = 0,$$

*where* $\varepsilon(d,X_j,Y_j)$ *and* $z(d,\mathbf{Y}^{(d)-},\mathbf{X}^{(d)-})$ *are as in* (10) *and* (12), *respectively.*

PROOF. First, let $T(x) = e^x \mathbf{1}_{(x<0)}$,

$$A(d, \mathbf{Y}^{(d)-}, \mathbf{X}^{(d)-}) = T\left(\sum_{j=1, j \neq i^*}^{d} \varepsilon(d, X_j, Y_j)\right) - T(z(d, \mathbf{Y}^{(d)-}, \mathbf{X}^{(d)-}))$$

and

$$\delta(d) = \left(\sum_{i=1}^{m} \mathrm{E}_{\mathbf{Y}_{\mathcal{J}(i,d)}^{(d)-}}[|W_i(d, \mathbf{X}_{\mathcal{J}(i,d)}^{(d)-}, \mathbf{Y}_{\mathcal{J}(i,d)}^{(d)-})|] + \sum_{i=1}^{m} c(\mathcal{J}(i,d))\ell^3 K \frac{d^{3\gamma_i/2}}{d^{3\alpha/2}}\right)^{1/2}.$$

We shall show that $A(d, \mathbf{Y}^{(d)-}, \mathbf{X}^{(d)-})|\mathbf{X}^{(d)-} \xrightarrow{p} 0$ and then use this result to prove convergence of expectations.

Similarly to the proof of Lemma A.7 in [10], we have

$$P_{\mathbf{Y}^{(d)-}}(|A(d, \mathbf{Y}^{(d)-}, \mathbf{X}^{(d)-})| \geq \delta(d))$$

(16)
$$\leq P_{\mathbf{Y}^{(d)-}}\left(\left|\sum_{j=1, j \neq i^*}^{d} \varepsilon(d, X_j, Y_j) - z(d, \mathbf{Y}^{(d)-}, \mathbf{X}^{(d)-})\right| \geq \delta(d)\right)$$

$$+ P_{\mathbf{Y}^{(d)-}}(-\delta(d) < z(d, \mathbf{Y}^{(d)-}, \mathbf{X}^{(d)-}) < \delta(d)).$$

By Markov's inequality and the proof of Lemma 8, the first term on the right-hand side is bounded by

$$\frac{1}{\delta(d)} \mathrm{E}_{\mathbf{Y}^{(d)-}}\left[\left|\sum_{j=1, j \neq i^*}^{d} \varepsilon(d, X_j, Y_j) - z(d, \mathbf{Y}^{(d)-}, \mathbf{X}^{(d)-})\right|\right] \leq \sqrt{\delta(d)} \xrightarrow{p} 0$$

as $d \to \infty$. Using conditioning and the proof of Lemma 9, the second term on the right-hand side becomes

$$P_{\mathbf{Y}^{(d)-}}(|z(d, \mathbf{Y}^{(d)-}, \mathbf{X}^{(d)-})| < \delta(d))$$

$$= \mathrm{E}_{\mathbf{Y}^{(n)-}}\left[\Phi\left(\frac{\delta(d) - \sum_{j=1, j \neq i^*}^{n} \varepsilon(d, X_j, Y_j) + \frac{\ell^2}{2} \sum_{i=1}^{m} R_i(d, \mathbf{X}_{\mathcal{J}(i,d)}^{(d)-})}{\ell \sqrt{\sum_{i=1}^{m} R_i(d, \mathbf{X}_{\mathcal{J}(i,d)}^{(d)-})}}\right)\right.$$

$$\left. - \Phi\left(\frac{-\delta(d) - \sum_{j=1, j \neq i^*}^{n} \varepsilon(d, X_j, Y_j) + \frac{\ell^2}{2} \sum_{i=1}^{m} R_i(d, \mathbf{X}_{\mathcal{J}(i,d)}^{(d)-})}{\ell \sqrt{\sum_{i=1}^{m} R_i(d, \mathbf{X}_{\mathcal{J}(i,d)}^{(d)-})}}\right)\right].$$

Using the convergence results developed in the proof of Lemma 9, along with the fact that $\delta(d) \xrightarrow{p} 0$ as $d \to \infty$ and the bounded convergence theorem, we deduce that the previous expression converges in probability to 0. Therefore, $A(d, \mathbf{Y}^{(d)-}, \mathbf{X}^{(d)-})|\mathbf{X}^{(d)-} \xrightarrow{p} 0$ and (15) follows by reapplying the bounded convergence theorem twice. □

6.4. *Simplified expression for the equivalent drift.*

LEMMA 11. *If condition* (5) *is satisfied, then*

$$\lim_{d \to \infty} \mathrm{E}_{\mathbf{X}^{(d)-}} \bigg[ \bigg| \mathrm{E}_{\mathbf{Y}^{(d)-}} [\exp\{z(d, \mathbf{Y}^{(d)-}, \mathbf{X}^{(d)-})\}; z(d, \mathbf{Y}^{(d)-}, \mathbf{X}^{(d)-}) < 0] \\ - \Phi\bigg(-\frac{\ell\sqrt{E_R}}{2}\bigg) \bigg| \bigg] = 0,$$

*where the functions* $\varepsilon(d, X_j, Y_j)$ *and* $z(d, \mathbf{Y}^{(d)-}, \mathbf{X}^{(d)-})$ *are as in* (10) *and* (12), *respectively.*

PROOF. The proof is similar to that of Lemma 9, the only difference lying in the fact that (Proposition 2.4 in [11])

$$\mathrm{E}_{\mathbf{Y}^{(d)-}}[\exp\{z(d, \mathbf{Y}^{(d)-}, \mathbf{X}^{(d)-})\}; z(d, \mathbf{Y}^{(d)-}, \mathbf{X}^{(d)-}) < 0] \\ = \mathrm{E}_{\mathbf{Y}^{(n)-}}\bigg[\exp\bigg(\sum_{j=1, j\neq i^*}^{n} \varepsilon(d, X_j, Y_j)\bigg) \\ \times \Phi\bigg(\frac{-\sum_{j=1,j\neq i^*}^{n} \varepsilon(d, X_j, Y_j) - \frac{\ell^2}{2}\sum_{i=1}^{m} R_i(d, \mathbf{X}^{(d)-}_{\mathcal{J}(i,d)})}{\sqrt{\ell^2 \sum_{i=1}^{m} R_i(d, \mathbf{X}^{(d)-}_{\mathcal{J}(i,d)})}}\bigg)\bigg].$$

□

**7. Discussion.** The theorems in this paper basically extend the i.i.d. work of Roberts, Gelman and Gilks [11] to a more general setting where the scaling term of each target component is allowed to depend on the dimension of the target distribution. The conclusions achieved are similar to those in [11], since the AOARs are identical; the sole difference lies in the optimal scaling values themselves. Condition (5), which says that no target component converges significantly faster than the others, ensures that the process behaves asymptotically as in the i.i.d. case. This work thus partially answers Open Problem #3 of [14].

These results can also be used to determine, for virtually any correlated multivariate normal target distribution, whether or not 0.234 is optimal. Contrary to what seemed to be a common belief, multivariate normal distributions do not always adopt a conventional limiting behavior and there exist cases where the AOAR is significantly smaller than 0.234 (see [1]).

It was shown in the i.i.d. case that although asymptotic, the results are fairly accurate in small dimensions ($d \geq 10$). In the present case, however, this fact is not always verified and care must be exercised in practice. In particular, if there exists a finite number of scaling terms such that $\lambda_j$ is close to $\alpha$ [but with $\lambda_j < \alpha$, otherwise condition (5) would be violated], then the optimal acceptance rate converges extremely slowly to 0.234 from above. For instance, suppose that $\boldsymbol{\Theta}^{-2}(d) = (d^{-\lambda}, 1, \ldots, 1)$ with $\lambda < 1$. The proposal scaling is then $\sigma^2(d) = \ell^2/d$ and the closer $\lambda$ is to 1, the slower the convergence of the optimal acceptance rate is to 0.234. In fact, for a multivariate normal target with $\lambda = 0.75$, simulations show that $d$ must be as large as 200,000 for the optimal acceptance rate to be reasonably close to 0.234; they also show that for $\alpha - \lambda \geq 0.5$, the asymptotic results are accurate in relatively small dimensions, just as in the i.i.d. case. Detailed examples and simulation studies illustrating the results introduced in this paper and in [1] are presented in [2].

## APPENDIX

LEMMA A.1. *Let $f$ be a $C^2$ probability density function (p.d.f.). If $(\log f(x))'$ is Lipschitz continuous, then $f'(x) \to 0$ as $x \to \pm\infty$.*

PROOF. The asymptotic behavior of a $C^2$ p.d.f. as $x \to \pm\infty$ can be one of three things: (1) $f(x) \to 0$, $f'(x) \to 0$; (2) $f(x) \to 0$, $f'(x) \not\to 0$; (3) $f(x) \not\to 0$, $f'(x) \not\to 0$. We prove that in cases (2) and (3), $(\log f(x))'$ is not Lipschitz continuous, which implies that (1) is the only possible option.

(2) $f(x) \to 0$, $f'(x) \not\to 0$: Since $f \to 0$, it follows that $\forall \epsilon > 0$, $\exists x_0(\epsilon) \in \mathbf{R}$ such that $\forall x \geq x_0(\epsilon)$, $f(x) < \epsilon$. Since $f' \not\to 0$, it follows that $\forall \epsilon > 0$, $\exists x^* \geq x_0(\epsilon) + 1$ such that $|f'(x^*)| > \limsup |f'|/2$. Because $f$ is $C^2$, we have $\forall 0 < \epsilon < \limsup |f'|/2$, $\exists y$ with $|x^* - y| \leq 1$ such that $|f'(y)| = \epsilon$. Now, choose $y^*$ to be the value of $y$ which minimizes $|x^* - y|$, but such that $f(y^*) > f(x^*)$. Given $0 < \epsilon < \limsup |f'|/2$, we then have

$$\sup_{x,y \in \mathbf{R}, x \neq y} \frac{|f'(x)/f(x) - f'(y)/f(y)|}{|x-y|} \geq \frac{||f'(x^*)|/f(x^*) - |f'(y^*)|/f(y^*)|}{1}$$
$$\geq \left| \frac{\limsup |f'|/2 - \epsilon}{\epsilon} \right|.$$

Since this is true for all $0 < \epsilon < \limsup |f'|/2$, the Lipschitz continuity assumption is violated.

(3) $f(x) \not\to 0$, $f'(x) \not\to 0$: Since $f$ is continuous, positive and $\int f = 1$, it follows that $\forall \epsilon > 0$, $\exists x_0(\epsilon) \in \mathbf{R}$ such that $f(x) < \epsilon$ for $x \geq x_0(\epsilon)$, except on a set $A_\epsilon$ of Lebesgue measure $\lambda(A_\epsilon) < \epsilon$. Since $(-\infty, \epsilon)$ is an open set, it follows that $B = \{x \in \mathbf{R} : f(x) < \epsilon\}$ must also be open; $A_\epsilon = B^c \cap [x_0(\epsilon), \infty)$ is then formed from closed intervals over which $f(x) \geq \epsilon$.

Since $f \not\to 0$, it follows that $\forall \epsilon > 0$, there exists an interval $[x(\epsilon), y(\epsilon)]$ in $A_\epsilon$ where the maximum value reached by $f$ over this interval ($h(\epsilon)$ say) is such that $h(\epsilon) > \limsup |f|/2$. There might be many values in the interval for which $f(x) = h(\epsilon)$, but all of these values will satisfy $f'(x) =$

0. Since $f(x(\epsilon)) = f(y(\epsilon)) = \epsilon$, it follows that $\sup_{x \in \mathbf{R}} f'(x) \geq \frac{h(\epsilon) - \epsilon}{y(\epsilon) - x(\epsilon)} > \frac{h(\epsilon) - \epsilon}{\epsilon}$. Hence, $\sup_{x \in \mathbf{R}} \frac{f'(x)}{f(x)} > \frac{h(\epsilon) - \epsilon}{\epsilon h(\epsilon)}$ and since this is true $\forall \epsilon > 0$, we have $\sup_{x \in \mathbf{R}} \frac{f'(x)}{f(x)} = \infty$. Given $\epsilon > 0$, we take $y$ to be one of the points in $[x(\epsilon), y(\epsilon)]$ such that $f(y) = h(\epsilon)$ and $f'(y) = 0$. We then have

$$\sup_{x,y \in \mathbf{R}, x \neq y} \frac{|f'(x)/f(x) - f'(y)/f(y)|}{|x - y|} \geq \sup_{x \in \mathbf{R}} \frac{|f'(x)/f(x) - 0|}{|x(\epsilon) - y(\epsilon)|}$$

$$> \sup_{x \in \mathbf{R}} \frac{|f'(x)/f(x) - 0|}{\epsilon} = \infty$$

and we see that the Lipschitz continuity assumption is violated. Note that in cases (2) and (3), we have considered the case where $x \to \infty$; we can construct a similar argument for the case where $x \to -\infty$. □

PROPOSITION A.2. *Let $\varepsilon(d, X_j, Y_j)$, $j = 1, \ldots, n$, be as in (10). If $\lambda_j < \alpha$, then $\varepsilon(d, X_j, Y_j) \underset{p}{\to} 0$.*

PROOF. By Taylor's theorem, we have for some $U_j \in (X_j, Y_j)$ or $(Y_j, X_j)$

$$E[|\varepsilon(d, X_j, Y_j)|]$$
$$= E[|(\log f(\theta_j(d) X_j))'(Y_j - X_j) + \tfrac{1}{2}(\log f(\theta_j(d) X_j))''(Y_j - X_j)^2$$
$$+ \tfrac{1}{6}(\log f(\theta_j(d) U_j))'''(Y_j - X_j)^3|].$$

Applying changes of variable and using the fact that $|(\log f(X))''|$ and $|(\log f(U))'''|$ are bounded by a constant, we obtain, for some $K > 0$,

$$E[|\varepsilon(d, X_j, Y_j)|] \leq \ell \frac{d^{\lambda_j/2}}{d^{\alpha/2}} K E[|(\log f(X))'|] + \left(\ell^2 \frac{d^{\lambda_j}}{d^\alpha} + \ell^3 \frac{d^{3\lambda_j/2}}{d^{3\alpha/2}}\right) K.$$

By assumption, $E[|(\log f(X))'|]$ is bounded by some finite constant. Since $\lambda_j < \alpha$, the previous expression converges to 0 as $d \to \infty$. To complete the proof of the proposition, we use Markov's inequality and find that for all $\epsilon > 0$, $P(|\varepsilon(d, X_j, Y_j)| \geq \epsilon) \leq E[|\varepsilon(d, X_j, Y_j)|]/\epsilon \to 0$ as $d \to \infty$. □

PROPOSITION A.3. *Let $R_i(d, \mathbf{X}_{\mathcal{J}(i,d)}^{(d)-})$ be as in (14), with $i \in \{1, \ldots, m\}$. We have $\sum_{i=1}^m R_i(d, \mathbf{X}_{\mathcal{J}(i,d)}^{(d)-}) \underset{p}{\to} E_R$, where $E_R$ is as in (6).*

PROOF. The expectation of each variable satisfies $E[R_i(d, \mathbf{X}_{\mathcal{J}(i,d)}^{(d)-})] = \frac{c(\mathcal{J}(i,d))}{d^\alpha} \frac{d^{\gamma_i}}{K_{n+i}} E[(\frac{f'(X)}{f(X)})^2]$. By independence between the $X_j$'s and the fact that $\text{Var}(X) \leq \text{E}[X^2]$, we obtain

$$\text{Var}\left(\sum_{i=1}^m R_i(d, \mathbf{X}_{\mathcal{J}(i,d)}^{(d)-})\right) \leq \sum_{i=1}^m \frac{1}{d^{2\alpha}} \frac{d^{2\gamma_i}}{K_{n+i}^2} c(\mathcal{J}(i,d)) \text{E}\left[\left(\frac{f'(X)}{f(X)}\right)^4\right].$$

By assumption, $\text{E}[(\frac{f'(X)}{f(X)})^4]$ is finite and since $c(\mathcal{J}(i,d))\, d^{2\gamma_i} < d^{2\alpha}$, the variance converges to 0 as $d \to \infty$. To conclude the proof, we use Chebychev's inequality and find that $\forall \epsilon > 0$, $\text{P}(|\sum_{i=1}^m R_i(d, \mathbf{X}_{\mathcal{J}(i,d)}^{(d)-}) - E_R| \geq \epsilon) \leq \frac{1}{\epsilon^2} \text{Var}(\sum_{i=1}^m R_i(d, \mathbf{X}_{\mathcal{J}(i,d)}^{(d)-})) \to 0$ as $d \to \infty$. $\square$

**Acknowledgments.** This work is part of my Ph.D. thesis; special thanks are due to my supervisor, Professor Jeffrey S. Rosenthal, without whom this work would not have been completed. His expertise, guidance and encouragement have been precious throughout my studies. I also acknowledge useful conversations with Professor Gareth O. Roberts and proofreading by Stephen B. Connor, as well as constructive comments from the anonymous referee which led to an improved paper.

DEPARTMENT OF STATISTICS
UNIVERSITY OF WARWICK
COVENTRY, CV4 7AL
UNITED KINGDOM
E-MAIL: M.Bedard@warwick.ac.uk